\documentclass[12pt]{article}
\usepackage{latexsym,amssymb,amsmath,amsfonts,amsthm,url,graphicx, caption}
\theoremstyle{definition}

\def\beq{ \begin{equation} }
\def\eeq{ \end{equation} }
\def\mn{\medskip\noindent}

\def\square{\vcenter{\vbox{\hrule height .4pt
  \hbox{\vrule width .4pt height 5pt \kern 5pt
        \vrule width .4pt} \hrule height .4pt}}}

\def\ans#1{}

\begin{document}

\title{Controlling the spread of \\
COVID-19 on college campuses}
\author{Molly Borowiak, Fayfay Ning, Justin Pei, Sarah Zhao,\\
 Hwai-Ray Tung, and Rick Durrett }

\maketitle

\begin{abstract}
This research was done during the DOMath program at Duke University from May 18 to July 10, 2020. At the time, Duke and other universities across the country were wrestling with the question of how to safely welcome students back to campus in the Fall. Because of this, our project focused on using mathematical models to evaluate strategies to suppress the spread of the virus on campus, specifically in dorms and in classrooms. For dorms, we show that giving students single rooms rather than double rooms can substantially reduce virus spread. For classrooms, we show that moving classes with size above some cutoff online can make the basic reproduction number $R_0<1$, preventing a wide spread epidemic. The cutoff will depend on the contagiousness of the disease in classrooms.
\end{abstract}

\section{Introduction}
As most readers know, the COVID-19 pandemic began in China in December 2019, then slowly spread around the world. Among its impacts on the US are changes to the education system. After Spring Break 2020, many colleges sent students home and switched to remote learning. During the summer, the big question became what to do in the fall. Colleges across the US have moved classes online, restricted the number of students on campus, and have modified semester schedules to run from mid-August to Thanksgiving with no breaks. To further control the spread of the virus on campus, they have made dining halls take-out only and reconfigured public spaces to encourage social distancing. In addition, there will be daily symptom reporting, increased testing, as well as quarantine and contact tracing for infected individuals. 

In this paper we will consider two methods to suppress the spread of
the virus. In Section \ref{sec:dorm}, we ask the question: How much
would spread be reduced if all students in dorms have single rooms? This
is a significant issue for Duke because undergraduates are required to
live on campus three of their four years, so at any time more than 80\%
live in dorms. By using a generalization of the household model
\cite{twolevels}, we show that the reduction can be substantial.

In Section \ref{sec:class}, we turn our attention to classrooms and ask
a two part question inspired by Gressman and Peck \cite{GrePeck}, who
analyzed a complex stochastic agent based model to determine the impact
of having large classes be online-only. The first part is: How does the
distribution of class sizes affect virus spread? To approach this
question, we considered two scenarios for a college that has $1000$ students, each
taking $3$ classes. In Scenario 1, all classes have size $30$. In
Scenario $2$, 1/4 of the classes have size 60, and 3/4 have size 20. By
calculating the maximum eigenvalue of the matrix that governs the spread
in the second case, we show that the basic reproduction number $R_0$ is 48\% larger in scenario 2 due to increased spread in larger
classes.

Knowing that larger classes promote virus spread, it is
natural to move classes above a certain threshold online. Thus, the
second question we ask in Section \ref{sec:class} is: Where should the
threshold be for moving classes online? To address this question, we
consider a scenario where 2405 students take three classes each and
there is one class of each size from $10$ through $120$. This is not a
very realistic situation but it does allow us to investigate a wide
range of cutoffs. We find that the threshold for moving classes online
to reduce the basic reproduction number $R_0 = 1$ and hence stop the
spread of the epidemic depends on the basic infection probability p in
the model. For the values we investigate, the threshold can be as large as 70 or as small as 30. 

\section{Household model and college dorms}
\label{sec:dorm}
The question that motivates this section is: To what extent would the spread of COVID-19 on the Duke campus be reduced if all students had single rooms? To answer this question, we will use results that have been developed for the household model. See \cite{twolevels} and references therein. The household model improves slightly on a homogeneous mixing model by partitioning the population into multiple households with the probability of infecting an individual being larger for individuals in the same household. 

We will review the necessary theory on the household model in Section \ref{sec:household_model}. In Section \ref{sec:single_rooms}, we examine the case where all students have single rooms by viewing dorms as households and applying the household model. In Section \ref{sec:double_rooms}, where we examine double rooms, we still view dorms as households, but generalize our approach on computing the dynamics within a household.

\subsection{Household model} \label{sec:household_model}

Let the population of size $N=mn$ be split into $m$ households of size $n$. For simplicity, we consider the case of fixed infection time. Let $p_H$ be the probability an individual $i$ infects individual $j$ who is another member of the same house. Let $p_G$ be the probability of infecting individual $j$ who is not in the same house as individual $i$. $H$ is used to denote household, and $G$ is used to denote global.

\mn
{\bf Step 1:} Connect each pair of individuals in the same house with probability $p_H$.

\mn
Let $\pi_k$ be the probability an individual belongs to a cluster of size $k$. Let $\mu = \sum_k k \pi_k$ be the mean cluster size, and let $G_H(z) =\sum_{k=1}^\infty \pi_k z^k$ be the generating function of the cluster size.

\mn
{\bf Step 2:} Make long distance connections between individuals in different households with probability $p_G$. 

\mn
The number of long distance infections caused by one person is $Y_1\approx$ Poisson($Np_G$) if $N$ is large. The generating function of $Y_1$ is $G_{\nu}(z) = \exp(Np_G(z-1))$.

From step 1 and 2, we see that
$$
R_0 = \mu N p_G.
$$
\noindent 
To compute the generating function, we recall that 
if $X,Y\ge 0$ are independent with generating functions $\phi_X$ and $\phi_Y$, then $X+Y$ has generating function $\phi_X\phi_Y$. This implies that the generating function of the number of individuals infected by a cluster of size $k$ is $\exp(Np_G(z-1))^k$ and hence that the number of individuals infected by the cluster of a randomly chosen individual is
$$
\sum_{k=0}^\infty \pi_k \exp(Np_G(z-1))^k = G_H(G_\nu(z)).
$$
A well-known result for branching processes implies that the probability of a large epidemic is the solution of
\begin{equation}
1- \zeta = G_H(G_\nu( \zeta)).
\label{eq:Zeta}
\end{equation}

When we apply the household model to our dorm models by viewing a dorm as a household, $G_H$ will vary between models since the choice of single or double rooms results in a slightly altered local spread within a dorm. On the other hand, since the global spread of the infection is unaffected by room size, $G_\nu$ will be the same for the two dorm scenarios.

\subsection{Single Rooms} \label{sec:single_rooms}

In adapting the household model to the spread of infection in a college dorm, $p_H$ will instead be written as $p_D$, since we are concerned with a dorm rather than a household. Similarly, $G_H$ will instead be written as $G_D^1$, with the superscript 1 indicating single rooms. 

We first compute $R_0$. Consider $m$ dorms with $n$ single rooms, where population size $N = mn$. If $n$ is large, e.g., 120 is the size of a typical freshman dorm at Duke, then the total number of individuals infected in one dorm as a result of one infected is the total progeny of a branching process. The branching process will have the offspring distribution Poisson($\lambda_1$), where $\lambda_1 = np_D$. We assume $np_D<1$ so that the mean of the total progeny is
\begin{equation}
\mu_1 = 1 + \lambda_1 + \lambda_1^2 + \cdots = \frac{1}{1-\lambda_1},
\label{eq:Mu_Single}
\end{equation}

\noindent 
Each individual infected in the dorm infects an average of $np_G$ individuals outside the dorm so the basic reproduction number is
\begin{equation}
R_0 = \mu_1 Np_G = \frac{N p_G}{1-np_D}.
\label{eq:R0_Single}
\end{equation}

Now, we solve for $G^1_D$. Breaking things down according to the number of individuals directly infected by the first individual, the generating function of the size of the epidemic within the dorm caused by one infected is 

\beq
G^1_D(z) = z\sum_{k=0}^{\infty} e^{-\lambda_1} \frac{\lambda_1^k}{k!}G^1_D(z)^k = z\exp(\lambda_1(G^1_D(z)-1)) = zG^1_b(G^1_D(z)),
\label{eq:G1D_notSolved}
\eeq

\noindent where $G^1_b(z) = \exp(np_D(z-1))$, the generating function of the spread of the infection within the dorm. One method for solving \eqref{eq:G1D_notSolved} is through iteration. We start with $h_0(z)=z$, which corresponds to 1 individual in generation 0,
and then iterate
\begin{equation}
h_{\ell}(z) = z\exp(\lambda_1(h_{\ell-1}(z)-1)).
\label{eq:Iter_Single}
\end{equation}

\noindent $h_{\ell}(z)$ is the generating function of the number of individuals infected in the first $\ell$ generations of the process, so $h_{\ell}(z) \downarrow h(z)$, a solution of the recursion.  

In the case of single rooms, an explicit solution for $G^1_D$ can be found by taking the derivative of equation \eqref{eq:G1D_notSolved} to get the differential equation 
$$
G^{1'}_D(z) = G^1_D(z)/z + \lambda_1 z G^1_D(z)G^{1'}_D(z),
$$
with the boundary condition $G^1_D(1) = 1$. It is well known that the solution to the differential equation is
\beq
G^1_D(z) = -W(zW^{-1}(-\lambda_1))/\lambda_1,
\label{lamW}
\eeq
where $W$ is the Lambert $W$ function which is defined by $W^{-1}(z) = ze^z$.

\begin{figure}[ht]
        \begin{center}
        \includegraphics[width=3.5in]{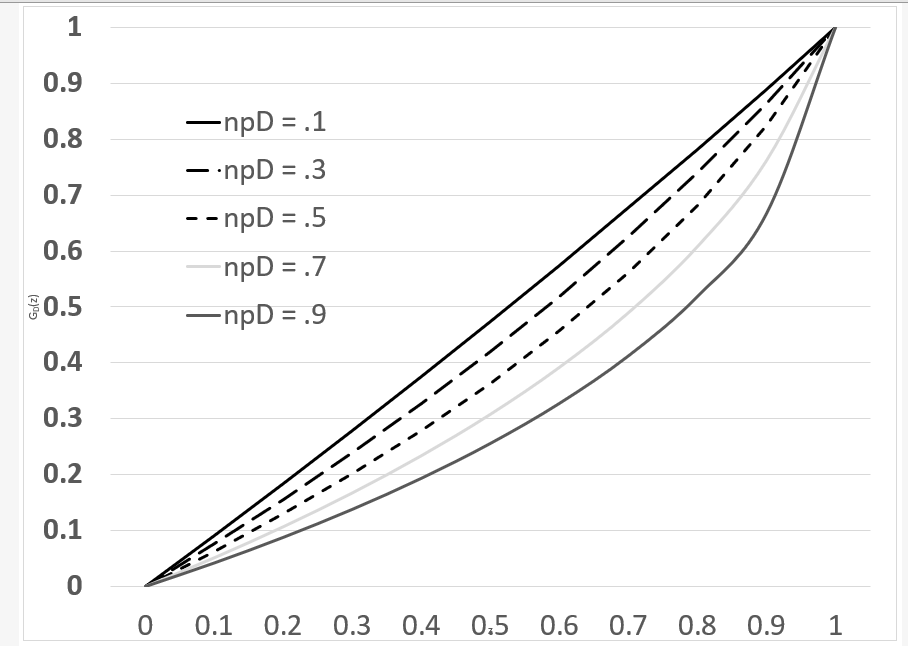}
        \caption{Generating function $G^1_D(z)$ for $np_D = .1, .3, .5, .7, .9$. Solved by iteration with Equation \eqref{lamW}. $z$ represents the initial proportion of students within the dorm that were infected.}
        \label{fig:GD Single}
        \end{center}
\end{figure}

\noindent
By \eqref{eq:Zeta}, the probability of a large epidemic is the largest solution of
\begin{equation}
1- \zeta = G^1_D(\exp(-Np_G\zeta)).
\label{eq:Zeta_Single} 
\end{equation}
This can be further simplified. Plugging in $G^1_D$ yields
$$
1-\zeta = -W(e^{-Np_G \zeta}(-\lambda_1 e^{-\lambda_1}))/\lambda_1 = -W(e^{-(Np_G + \lambda_1)\zeta}e^{\lambda_1(\zeta-1)}(-\lambda_1))/\lambda_1.
$$
Multiplying through by $-\lambda_1$, then applying $W^{-1}$ to both sides gives
$$
\lambda_1(\zeta-1)e^{\lambda_1(\zeta-1)} = e^{-(Np_G + \lambda_1)\zeta}e^{\lambda_1(\zeta-1)}(-\lambda_1).
$$
Recalling that $\lambda_1 = np_D$ and simplifying, we obtain 
$$
(1-\zeta) = e^{-(Np_G +np_D)\zeta},
$$
which is the survival probability for a branching process with a Poisson($R_0 = Np_G +np_D$) distribution. We do not have an intuitive understanding of why this is true.

\subsection{Double Rooms} \label{sec:double_rooms}
Suppose now that we have $m_1$ dorms with $n_1$ double rooms, so that $N = 2m_1n_1$. For the sake of comparison, we let the number of people in a dorm be the same, regardless of whether they have only single rooms or only double rooms, i.e., $n = 2n_1$. Let $p_L$ be the probability of infecting your roommate, and $p_D$ be the probability of infecting an individual within your dorm. To find the analogue of $G^1_b$ in this scenario, we think of a generation as having a first step where the ancestor splits into two with probability $p_L$, stays 1 with probability  $1-p_L$ and then a second step in which long range connections are formed within the dorm. The offspring distribution is now a mixture of two Poissons:
$$
p_L \hbox{Poisson}(2 \cdot 2n_1p_D)  + (1-p_L)\hbox{Poisson}(2n_1 p_D).
$$
The mean number of infections in one generation, $\lambda_2$, will be the sum of the means of the two Poisson distributions so $\lambda_2 = (1+p_L)2n_1p_D$.

Imitating calculations in Section \ref{sec:single_rooms}, if $\lambda_2 < 1$, the mean of the total progeny of the branching process is
\begin{equation}
\mu_2 = 1 + \lambda_2 + \lambda_2^2 + \cdots = \frac{1}{1-\lambda_2},
\label{eq:Mu_Double}
\end{equation}

\noindent and the basic reproduction number is
\begin{equation}
R_0 = \mu_2 Np_G = \frac{Np_G}{1 - (1+p_L)2n_1p_D}.
\label{eq:R0_Double}
\end{equation}
To compute $G^2_D$, the analogue of $G^1_D$ for the double room scenario, we note that the generating function of the number of infections in one generation is
$$
G^2_b(z) = p_L \exp(2n_1p_D(z-1))^2 + (1-p_L) \exp(2n_1p_D(z-1)), 
$$

\noindent so $G^2_D$ satisfies
\begin{equation}
G^2_D(z) = zG^2_b(G^2_D(z)).
\label{eq:G2D_notsolved}
\end{equation}

\noindent $G^2_D$ can then be solved with the iteration method described in Section \ref{sec:single_rooms}.

\begin{figure}[ht]
    \begin{center}
    \includegraphics[width=3.5in]{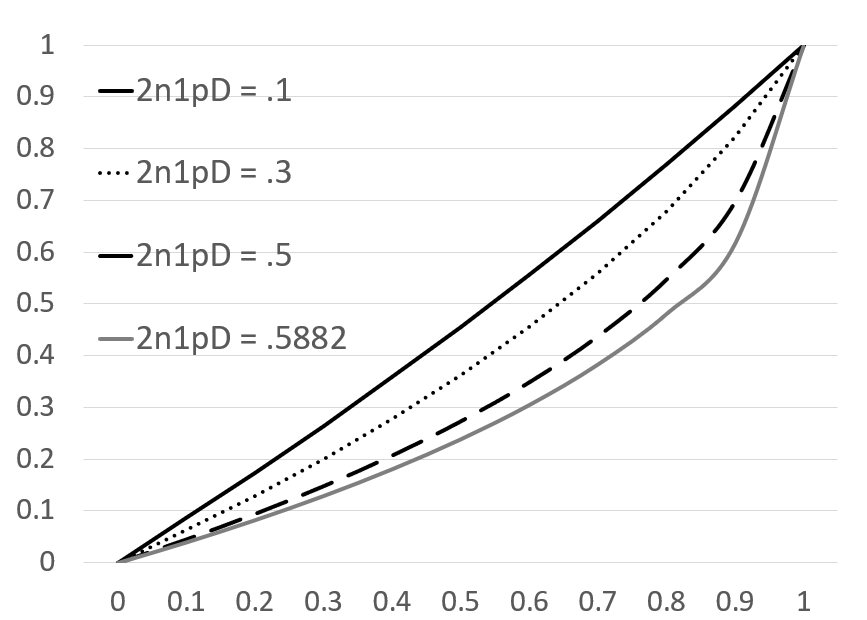}
    \caption{Generating function $G^2_D(z)$ for $2n_1p_D = .1, .3, .5, .5882$ while $p_L = .7$. Solved by iteration. When $p_L = .7$ and $\lambda_2 > .5882$, there is a large epidemic within the dorm. The graph of $G^2_D(z)$ increased with $z$, regardless of the value of $2n_1p_D$. As $2n_1p_D$ increased, the value of $G^2_D(z)$ at respective values of z slightly decreased.}
    \label{fig:GD Double}
    \end{center}
\end{figure}

Having computed $G^2_D$ we can compute the probability of a large epidemic as in the previous case by numerically solving
\begin{equation}
1- \zeta = G^2_D(\exp(-Np_G\zeta)).
\label{eq:Zeta_Double}
\end{equation}

\begin{figure}[ht]
\centering
\includegraphics[width=65mm]{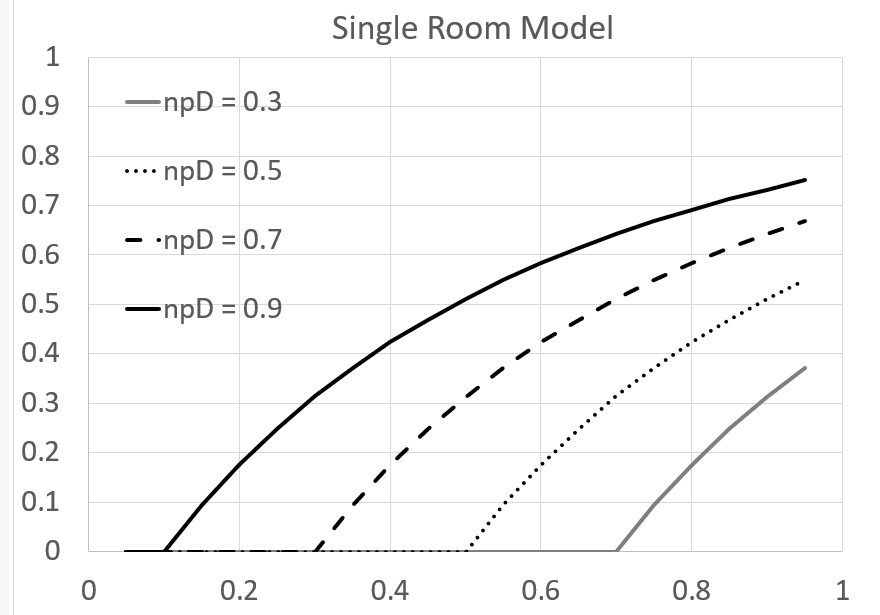}
\hspace{5mm}
\includegraphics[width=65mm]{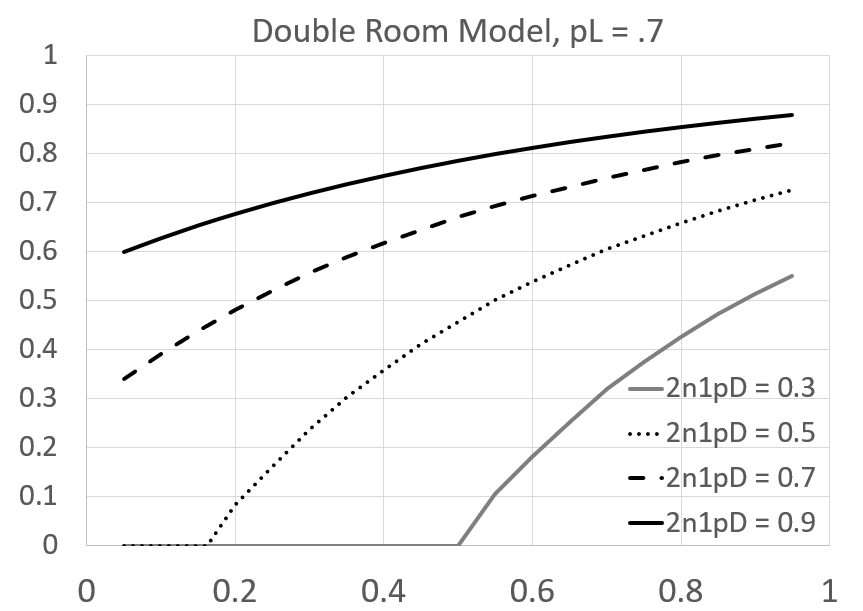}
\caption{Size of the epidemics in the single and double room scenarios of the household model when $np_D = 2n_1p_D = 0.3, 0.5, 0.7, 0.9$ and $Np_G$ varies.  
}
\label{fig:Zeta}
\end{figure}

A comparison of the epidemic probabilities $\zeta$ for  double-room and single-room dorms is given in Figure \ref{fig:Zeta}. Note that in the double-room situation when $2n_1p_D = 0.7$ or 0.9, there is a positive probability of a large epidemic within a single dorm, so the curves are positive when $Np_G=0$. When $2n_1p_D = 0.3$ or 0.5, $p_G$ needs to exceed a positive threshold for there to be an epidemic in a university with double-room dorms, but those thresholds are higher when there are only single-room dorm.
	      
\clearpage
\section{Spread of infection in college classes} \label{sec:class}

This section addresses two questions. The first asks how different class size distributions affect pandemic spread. The second asks whether taking classes above a certain size online can prevent outbreaks. 

To address these questions, we consider some simplified models in which each student takes three classes. The situation can be described by a graph in which there are $n$ students and $m$ classes with sizes $c_i$ $1\le i \le m$, and there is an edge from each student to the three classes they are enrolled in. Each classroom is assumed to be homogeneously mixing. To formulate the dynamics, we say that the three class members that correspond to one student are their alter egos.

An individual that is infected in class $i$ at time, adds his two alter egos (his presence in other classes) to the infected population at time $1/2$ and then infections at time $1$ are produced. None of the infecteds from time $1/2$ are present at time $1$.

\subsection{Equal class sizes reduce spread} 
\label{sec:scen1and2}
To show that class sizes impact spread, we compare two scenarios for a college with $1000$ that each take three classes. 

\mn
{\bf Scenario 1.} 100 classes of size 30. An individual infected in class $1$ at time 0 will produce a Poisson($29p$) number  of new infected individuals in class $1$ at time 1. His two alter egos are equally likely to be in classes 2-100, and each alter ego will produce Poisson($29p$) infecteds in whichever class they are in. So if $m_{i,j}$ is the mean number of infecteds in class $j$ at time 1 from a single infected person in class $i$, then 
$$
m_{i, i} = 29 p \qquad m_{i,j} = \frac{2}{99} \cdot 29p \quad\hbox{when $j \neq i$}
$$
From this we see that an infected at time 0 produces an average of $3\cdot 29p$ infecteds at time 1 so
\beq
R_0 = 87p.
\label{R01}
\eeq

\mn
{\bf Scenario 2}. 25 classes of size 60 ($1 \le i \le 25$) and 75 of size 20  ($26 \le i \le 100$). 
Let $m(i,j)$ be the expected number of infections in class $j$ caused by one infected student in $i$. $c_i$ and $c_j$ represent that sizes of class i and j respectively. In this scenario, the probability that an individual in class $i$ has an alter ego in a different class $j$ is $2c_j/(3000-c_i)$, where $c_i$ is the size of class $i$. The $2$ comes from $2$ alter egos. Thus, we have 
\begin{align*}
c_i = \begin{cases}
60 &  i \leq 25\\
20 &  i \geq 26
\end{cases}, \qquad
m_{i, j}&  = \begin{cases}
(c_i - 1)p &  i = j\\
2 \cdot \frac{c_j}{3000 - c_i} \cdot (c_j - 1)p &  i \neq j
\end{cases} 
\end{align*}
\noindent Since the eigenvector $x$ associated with the maximum eigenvalue will
have $x_i=a$ for $1\le i \le 25$ and $x_i=b$ for $26\le i \le 100$,
the eigenvalue problem for the $100 \times 100$ matrix can be reduced to
the eigenvalue problem for the $2\times 2$ matrix with entries
$$
\begin{bmatrix}
59 + 24 \cdot 59 \cdot 2 \cdot 60/2940 & 75 \cdot 19 \cdot 2 \cdot 20/2940\\
25 \cdot 59 \cdot 2 \cdot 60/2980 & 19 + 74 \cdot  19 \cdot 2 \cdot 20/2980
\end{bmatrix}p
= \begin{bmatrix}
116.796 & 19.388\\
59.396 & 37.783
\end{bmatrix}p
$$

\noindent Solving a quadratic shows that this matrix has eigenvalues
$129.380p$ and $25.288p$. Thus,
\beq
R_0 = 129.38p.
\label{R02}
\eeq
\noindent This value is 48\% larger than the $R_0$ for scenario 1, and suggests that larger classes increase virus transmission. 

\subsection{Stochastic model for scenarios 1 and 2}

For a more detailed look at the two scenarios, we turn to a stochastic model. As in the previous section the model takes place in discrete time. Each time step is roughly one class meeting, so that two time steps are roughly equivalent to a week in real time. An infected person in a class infects each other individual with probability $p-0.01$ and has a 0.5 probability of entering quarantine. Using \eqref{R01} we see that $R_0=2(0.87) = 1.74$ in scenario 1. \eqref{R02} implies that $R_0=2(1.2938) = 2.5876$. $R_0=2.5$ is a commonly used value for an uncontrolled COVID-19 epidemic.  Figures 4 and 5 confirm that Scenario 2 is worse than Scenario 1.

\begin{figure}[ht]
\begin{center}
\includegraphics[width=75mm]{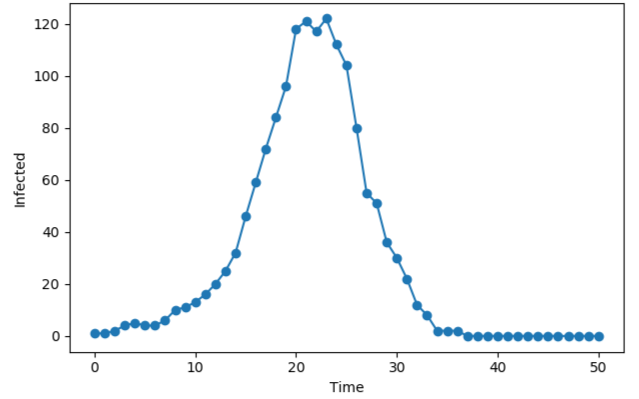}
\includegraphics[width=75mm]{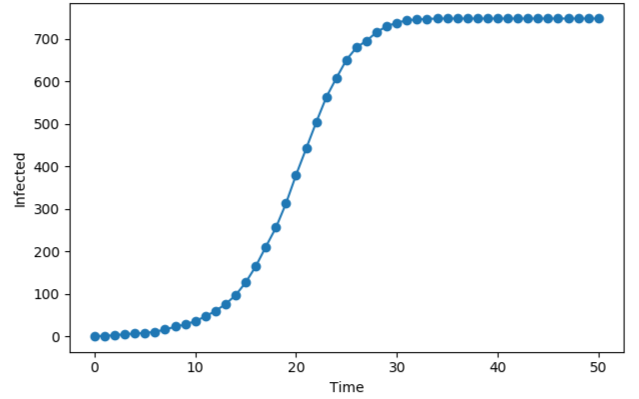}
\caption{A typical realization of the epidemic when all classes have size 30. The left panel gives the current number of infections. The right the cumulative number. In this situation, it takes roughly 11 weeks for the epidemic to reach its peak. By week 15, almost all individuals on campus have been infected}
\label{fig:exsame}
\end{center}
\end{figure}

The next figure compares the equal class size and 60--20 scenarios. In both scenarios, the outbreak generally grows out of control within a few weeks and infects the majority of students. In Scenario 2, the outbreaks tend to develop faster, have a higher peak, and vary less in size and dynamics. Scenario 2 also produces a higher probability a large outbreak will occur from one individual, and larger final epidemic sizes. Thus, it is clear that the presence of larger classes, even if average class size remains the same, increases the risks for students on campus.

\clearpage 

\begin{figure}[ht]
\centering
\includegraphics[width=75mm]{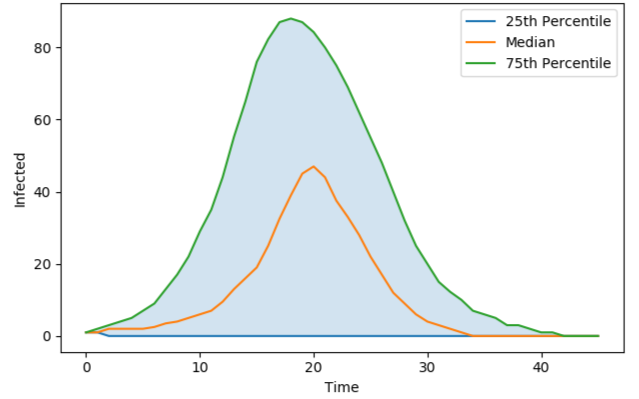}
\includegraphics[width=75mm]{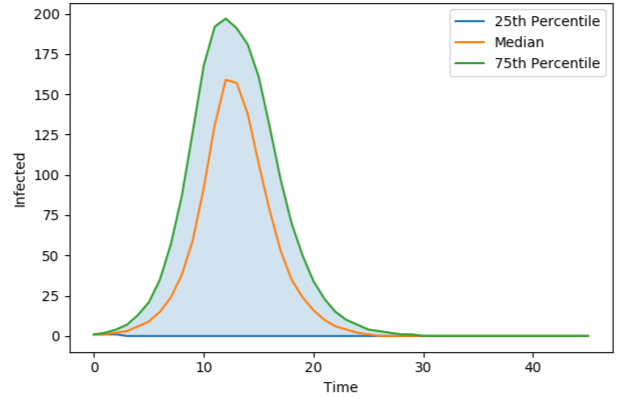}
\caption{Sample epidemics generated by Python. We ran the model $1,000$ times for each scenario. The left panel gives the result for scenario 1, the right for scenario 2. The top line is the 75th percentile, the next line is the median, the one parallel to the axis is the 25th percentile. Note that, in Scenario 1, the probability of no community infection is $0.243$, and the probability of no major outbreak (no more than 20 total infected) is $0.432$. For Scenario 2, the probabilities are $0.203$ and $0.309$ respectively.}
\label{fig:sim_time_plot}
\end{figure}

\begin{figure}[ht]
\centering
\includegraphics[width=0.5\linewidth]{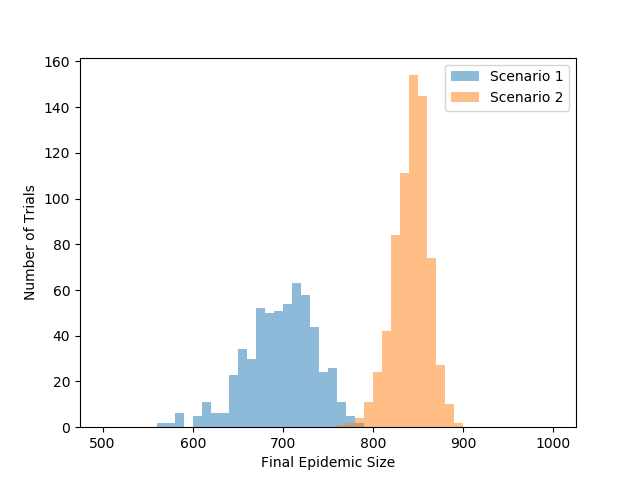}
\caption{Distributions of final epidemic size given that a large outbreak occurs. The data is the same data shown in Figure \ref{fig:sim_time_plot}.}
\end{figure}

\subsection{Moving Classes Online} 
In this section, we consider the effects of moving classes above a given size online. To do so, we consider a university with 111 classes, one each of sizes ranging from 10, 11... 120. The total class size is $10+11+...+120 = 7215$ and there are $7215/3 = 2405$ students. To begin, we first consider when all classes are in person.

Following the calculations for scenario 2 in Section \ref{sec:scen1and2}, we see that the average number fo people infected in classroom $j$ by an infected individual in classroom $i$ is
            \begin{align*}
                m_{i,j} = \left\{ \begin{array}{cc} 
                                (c_j - 1) p & \hspace{5mm} i=j \\ 
                                2 \cdot \frac{c_j}{7215 -c_i} \cdot (c_j-1) p & \hspace{5mm} i \neq j ,\\
                                \end{array} \right.
            \end{align*}
Similar to the previous scenarios in this section, we find the largest eigenvalue of the above matrix to determine $R_0$. Although we were unable to find the eigenvalue analytically, we used Python to determine that 

$$R_0 = 251.5p.$$

\noindent 
When $p = 0.01$, $R_0 = 2.515$.

\begin{figure}[ht]
	\begin{center}
    \includegraphics[width=4in]{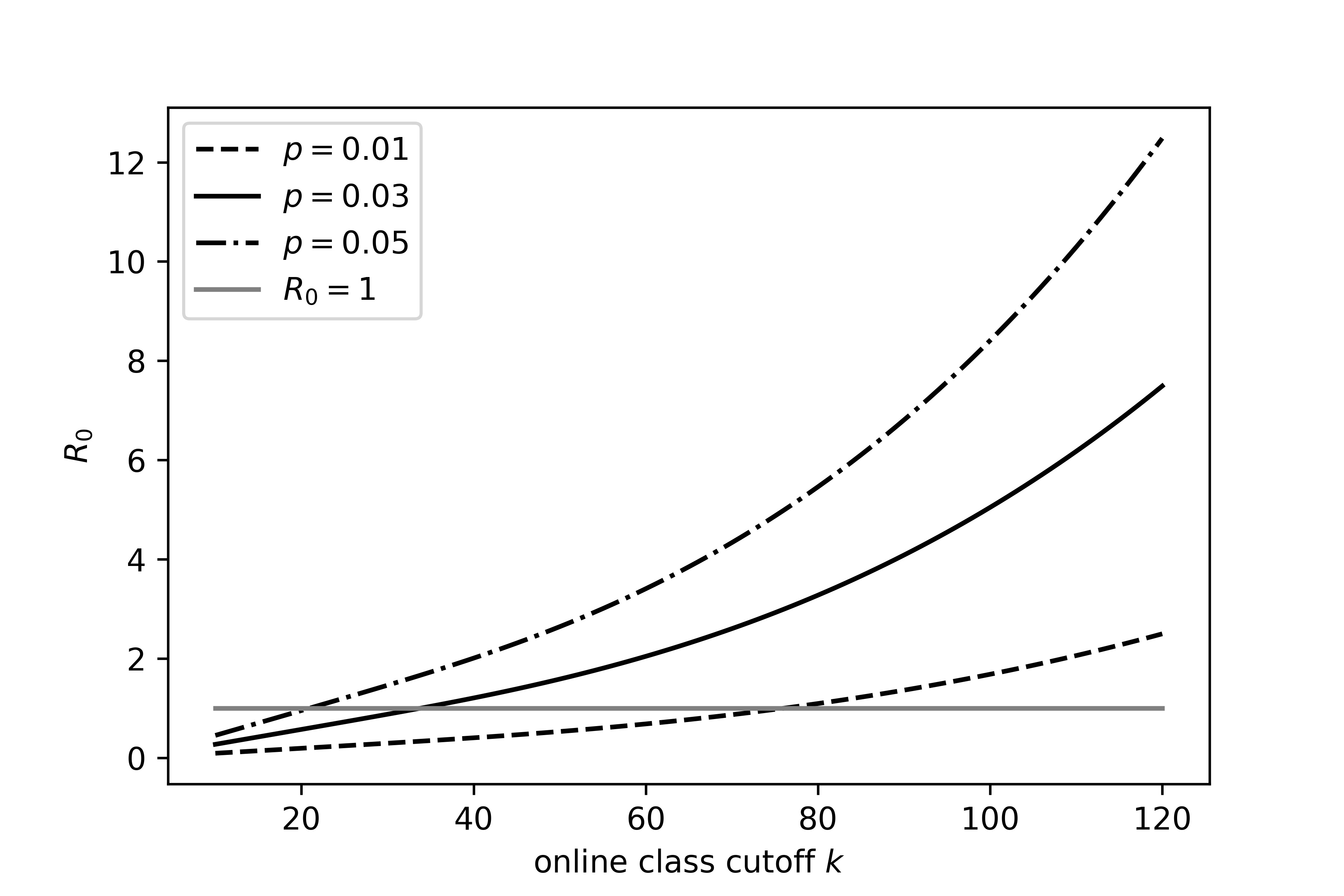}
    \caption{Basic reproduction number of COVID-19 after moving classes of size $>k$ online for varying values of $p$. Establishing a sufficiently small cutoff for moving classes online can make $R_0 \leq 1$, preventing an outbreak for a large range of $p$. What is considered sufficiently small depends greatly on $p$. }
    \label{fig:pyimage}
    \end{center}
\end{figure}

We now look at the situation where the university moves all of its classes of size $> k$. Viruses can't spread through online classes, we have
$$ m_{i,j} = 0 \text{ when } c_j > k.$$
The other values of $m(i, j)$ remain unchanged. The change in $R_0$ depending on $k$ and various $p$ can be seen in Figure \ref{fig:pyimage}, which makes clear that establishing a cutoff can sharply reduce virus spread through classes. Although Gressman and Peck suggested a cutoff of $30$ based on their stochastic simulations for spread in universities \cite{GrePeck}, our results in Figure \ref{fig:pyimage} suggest that the optimal cutoff depends largely on $p$. For sufficiently small $p$, even cutoffs as high as $70$ can be effective. 

\end{document}